\documentclass[11pt,fullpage]{article}

\newcommand{\ol}{\setlength{\itemsep}{0pt.}\begin{enumerate}}
\newcommand{\eol}{\end{enumerate}\setlength{\itemsep}{-\parsep}}
\newcommand{\ignore}[1]{}
\def\R{{\bf R}}

\title{\bf A short, based on the mixed volume, proof of Liggett's theorem on the convolution
of ultra-logconcave sequences}
 
\author{Leonid Gurvits \thanks{%
{\tt gurvits@lanl.gov}. Los Alamos National Laboratory, 
Los Alamos, NM. } 
}

\begin{document}



\maketitle

\begin{abstract}
R. Pemantle conjectured, and T.M. Liggett proved in 1997, that
the convolution of two ultra-logconcave is ultra-logconcave.
Liggett's proof is elementary but long. We present here
a short proof, based on the mixed volume of convex sets.

\end{abstract} 

 
\newtheorem{THEOREM}{Theorem}[section]
\newenvironment{theorem}{\begin{THEOREM} \hspace{-.85em} {\bf :} 
}%
                        {\end{THEOREM}}
\newtheorem{LEMMA}[THEOREM]{Lemma}
\newenvironment{lemma}{\begin{LEMMA} \hspace{-.85em} {\bf :} }%
                      {\end{LEMMA}}
\newtheorem{COROLLARY}[THEOREM]{Corollary}
\newenvironment{corollary}{\begin{COROLLARY} \hspace{-.85em} {\bf 
:} }%
                          {\end{COROLLARY}}
\newtheorem{PROPOSITION}[THEOREM]{Proposition}
\newenvironment{proposition}{\begin{PROPOSITION} \hspace{-.85em} 
{\bf :} }%
                            {\end{PROPOSITION}}
\newtheorem{DEFINITION}[THEOREM]{Definition}
\newenvironment{definition}{\begin{DEFINITION} \hspace{-.85em} {\bf 
:} \rm}%
                            {\end{DEFINITION}}
\newtheorem{EXAMPLE}[THEOREM]{Example}
\newenvironment{example}{\begin{EXAMPLE} \hspace{-.85em} {\bf :} 
\rm}%
                            {\end{EXAMPLE}}
\newtheorem{CONJECTURE}[THEOREM]{Conjecture}
\newenvironment{conjecture}{\begin{CONJECTURE} \hspace{-.85em} 
{\bf :} \rm}%
                            {\end{CONJECTURE}}
\newtheorem{PROBLEM}[THEOREM]{Problem}
\newenvironment{problem}{\begin{PROBLEM} \hspace{-.85em} {\bf :} 
\rm}%
                            {\end{PROBLEM}}
\newtheorem{QUESTION}[THEOREM]{Question}
\newenvironment{question}{\begin{QUESTION} \hspace{-.85em} {\bf :} 
\rm}%
                            {\end{QUESTION}}
\newtheorem{REMARK}[THEOREM]{Remark}
\newenvironment{remark}{\begin{REMARK} \hspace{-.85em} {\bf :} 
\rm}%
                            {\end{REMARK}}
\newtheorem{FACT}[THEOREM]{Fact}
\newenvironment{fact}{\begin{FACT} \hspace{-.85em} {\bf :} 
\rm}%
		            {\end{FACT}}

 
\newcommand{\thm}{\begin{theorem}}
\newcommand{\lem}{\begin{lemma}}
\newcommand{\pro}{\begin{proposition}}
\newcommand{\dfn}{\begin{definition}}
\newcommand{\rem}{\begin{remark}}
\newcommand{\xam}{\begin{example}}
\newcommand{\cnj}{\begin{conjecture}}
\newcommand{\prb}{\begin{problem}}
\newcommand{\que}{\begin{question}}
\newcommand{\cor}{\begin{corollary}}
\newcommand{\fac}{\begin{fact}}

\newcommand{\prf}{\noindent{\bf Proof:} }
\newcommand{\ethm}{\end{theorem}}
\newcommand{\elem}{\end{lemma}}
\newcommand{\epro}{\end{proposition}}
\newcommand{\edfn}{\bbox\end{definition}}
\newcommand{\erem}{\bbox\end{remark}}
\newcommand{\exam}{\bbox\end{example}}
\newcommand{\ecnj}{\bbox\end{conjecture}}
\newcommand{\eprb}{\bbox\end{problem}}
\newcommand{\eque}{\bbox\end{question}}
\newcommand{\ecor}{\end{corollary}}
\newcommand{\efac}{\end{fact}}
\newcommand{\eprf}{\bbox}
\newcommand{\beqn}{\begin{equation}}
\newcommand{\eeqn}{\end{equation}}
\newcommand{\wbox}{\mbox{$\sqcap$\llap{$\sqcup$}}}
\newcommand{\bbox}{\vrule height7pt width4pt depth1pt}
\newcommand{\qed}{\bbox}

\newcommand{\rarrow}{\rightarrow}
\newcommand{\larrow}{\leftarrow}
\newcommand{\grad}{\bigtriangledown}

\overfullrule=0pt
\def\setof#1{\lbrace #1 \rbrace}
\section{Introduction}
Let ${\bf a} = (a_0,...,a_m)$ and ${\bf b} = (b_0,...,b_n)$ be two real sequences.
Their convolution ${\bf c} ={\bf a} \star {\bf b}$ is defined as
$c_k = \sum_{i+j = k} a_i b_j, 0 \leq k \leq n+m$. 
A nonnegative sequence ${\bf a} = (a_0,...,a_m)$ is said to be {\it  logconcave} if
\beqn \label{lc}
a_i^{2} \geq a_{i-1} a_{i+1}, 1 \leq i \leq m-1.
\eeqn

Following Permantle and \cite{lig},
we say that a nonnegative sequence ${\bf a} = (a_0,...,a_m)$ is 
{\it ultra-logconcave of order} $d \geq m$ ($ULC(d)$) if the sequence
$\frac{a_i}{{d \choose i}}, 0 \leq i \leq m$ is {\it  logconcave}, i.e. 
\beqn \label{ulc} 
(\frac{a_i}{{d \choose i}})^2 \geq \frac{a_{i-1}}{{d \choose i-1}} \frac{a_{i+1}}{{d \choose i+1}}, 1 \leq i \leq m-1.
\eeqn
The next result was conjectured by R. Pemantle and proved by T.M. Liggett in 1997 \cite{lig}.

\thm \label{tom}
The convolution of a $ULC(l)$ sequence ${\bf a}$ and a $ULC(d)$ sequence ${\bf b}$ is $ULC(l+d)$.
\ethm

\rem \label{pert}
It is easy to see, by a standard pertubration argument, that it is sufficient to consider a positive case:
$$
{\bf a} =(a_0,...,a_l); a_i > 0, 0 \leq i \leq l \quad \mbox{and} \quad {\bf b} =(b_0,...,b_d); b_i > 0, 0 \leq i \leq d.
$$
The (simple) fact that the convolution of {\it  logconcave} sequences is also {\it  logconcave} was proved in \cite{Polya}
in 1949.
\erem

We present in this paper a short proof of Theorem(\ref{tom}).

\section{The mixed volume}
Let ${\bf K} = (K_1...K_n)$ be a $n$-tuple of convex compact subsets in the Euclidean space $\R^n$, and
let $V_{n}(\cdot)$ be the Euclidean volume in $\R^n$. It is well known Herman Minkowski result
(see for instance \cite{BZ88}), that the value of the  
$V_{n}(\lambda_1 K_1 + \cdots
\lambda_n K_n)$ is a homogeneous polynomial of degree $n$, called the Minkowski polynomial, in
nonnegative variables $\lambda_1...\lambda_n$, where $''+''$
denotes Minkowski sum, and $\lambda K$ denotes the dilatation of
$K$ with coefficient $\lambda$.  
The coefficient $V(K_1...K_n)$ of
$\lambda_1\cdot \lambda_2\ldots \cdot \lambda_n$ is called the {\it mixed
volume} of $K_1...K_n$. Alternatively,
$$
V(K_1...K_n) = \frac{\partial^n}{\partial \lambda_1...\partial
\lambda_n} V_{n}(\lambda_1 K_1 + \cdots
\lambda_n K_n).
$$
{\it The Alexandrov-Fenchel inequalities state that}
\beqn
V(K_1,K_2,...,K_n)^2 \geq V(K_1,K_1,...,K_n)V(K_2,K_2,...,K_n)
\eeqn

It follows that if $P,Q \subset R^n$ are convex compact sets then
$$
Vol_{n}(t P + Q) = \sum_{0 \leq i \leq n} a_i t^i,
$$
where 
$$
a_0 = Vol_{n}(Q) = \frac{1}{n!}V(Q,\cdots,Q), a_1 = \frac{1}{(n-1)! 1!} V(P,Q, \cdots,Q),..., a_n = Vol_{n}(Q) = \frac{1}{n!} V(P,\cdots,P)
$$
Which means that the sequence $(a_0,...,a_n)$ is $ULC(n)$.\\
The next remarkable result was proved by G.S. Shephard in 1960:
\thm \label{shep}
A sequence $(a_0,...,a_n)$ is $ULC(n)$ if and only if there exist two convex compact  sets $P,Q \subset R^n$ such that
$$
\sum_{0 \leq i \leq n} a_i t^i = Vol_{n}(t P + Q), t \geq 0
$$
\ethm

\rem
The "if" part in Theorem(\ref{shep}), which is a particular case of {\it the Alexandrov-Fenchel inequalities}, is not simple,
but was proved seventy years ago \cite{Al38}. The proof of the "only if" part in Theorem(\ref{shep}) is not difficult and short.
G.S. Shephard first considers the case of positive coefficients, which is already sufficient for our application.
In this positive case one chooses $Q = \{(x_1,...,x_n) : \sum_{1 \leq i \leq n} x_i \leq 1; x_i \geq 0\}$.
In other words, the set $Q$ is the standard simplex in $R^n$. And the convex compact set 
$P = Diag(\lambda_1,...,\lambda_n) Q, \lambda_1 \geq...\geq \lambda_n > 0$.
The general nonnegative case case handled by the topological theory of convex compact subsets. 
\erem

\section{Our proof of Theorem(\ref{tom})}
\prf
Let ${\bf a} = (a_0,...,a_l)$ be $ULC(l)$ and ${\bf b} =(b_0,...,b_d)$ be $ULC(d)$. Define two univariate
polynomials $R_1(t) = \sum_{0 \leq i \leq l} a_i t^i$ and $R_2(t) = \sum_{0 \leq j \leq c} a_i t^j$.\\
Then the polynomial $R_1(t) R_2(t) := R_3(t) = \sum_{0 \leq k \leq l+d} c_k t^k$, where the
sequence\\
${\bf c} = (c_0,...,c_{l+d})$ is the convolution, ${\bf c} = {\bf a} \star {\bf b}$.\\
It follows from the "only if" part of Theorem(\ref{shep}) that
$$
R_1(t) = Vol_{l} (t K_1 + K_2) \quad \mbox{and} \quad  R_2(t) = Vol_{d}(t C_1 + C_2),
$$
where $K_1, K_2, C_1, C_2$ are convex compact sets; $K_1, K_2 \subset R^l$ and $C_1, C_2 \subset R^d$.\\
Define the next two convex compact subsets of $R^{l+d}$:
$$
P = K_1 \times C_1 \quad \mbox{and} \quad Q = K_2 \times C_2.
$$
Here the cartesian product $A \times B$ of two subsets $A \subset R^l$ and $B \subset R^d$ is defined as
$$
A \times B := \{ (X,Y) \in R^{l+d} : X \in A, Y \in B \}.
$$
Clearly, the Minkowski sum $tP + Q = (tK_1 + K_2) \times (tC_1 + C_2), t \geq 0$.\\
It follows that $Vol_{l} (t K_1 + K_2) Vol_{d}(t C_1 + C_2) = Vol_{l+d}(t P +Q)$. Therefore
the polynomial\\ 
$R_{3}(t) = Vol_{l+d}(t P +Q)$. It follows from {\it the Alexandrov-Fenchel inequalities} (the "if" part
of Theorem(\ref{shep}))
that the sequence of its coefficients ${\bf c} = {\bf a} \star {\bf b}$ is $ULC(l+d)$.
\eprf
\section{Final comments}
\begin{enumerate}
\item Let ${\bf a} = (a_0,...,a_m)$ be a real sequence, satisfying the Newton inequalities of order $m$ (\ref{ulc}). I.e. we dropped
the condition of nonnegativity from the definition of {\it ultra-logconcavity}. It is not true that ${\bf c} = {\bf a} \star {\bf a}$ satisfies the Newton inequalities of order $2m$.\\
Indeed, consider ${\bf a} = (1,a,0,-b,1)$, where $a,b > 0$ and $\frac{b^2}{a^2}$ is sufficiently small.\\
Then $c_6 = b^2, c_5 = 2 a, c_4 = 2(1 - ab)$ and the number 
$$
\frac{c_{5}^{2}}{c_4 c_6} = 2 \frac{a^2}{b^2 (1 - ab)}
$$
converges to zero if the numbers $a,b, \frac{a}{b}$ converge to zero. 
\item
The reader can found further implications and generalization of Theorem(\ref{shep}) in \cite{mtns}.
\end{enumerate}

\end{document}